\documentclass{amsart}
\usepackage{amscd, amsmath, amsthm, amssymb}

\usepackage[bookmarksnumbered,colorlinks, plainpages]{hyperref}
\usepackage{stackrel}
\hypersetup{
pdftitle={Research Paper},
pdfauthor={Reza Abdolmaleki},
pdfsubject={Defining ideals of Cohen-Macaulay fiber cones},
pdfcreator={Reza Abdolmaleki},
colorlinks,
linkcolor=blue,
citecolor=magenta,
anchorcolor=red,
bookmarksopen,
urlcolor=black,
filecolor=red,
}
\theoremstyle{plain}
\newtheorem{Theorem}{Theorem}[section]

\newtheorem{thm}[Theorem]{Theorem}
\newtheorem{lem}[Theorem]{Lemma}

\newtheorem{prop}[Theorem]{Proposition}

\theoremstyle{definition}

\newtheorem{ex}[Theorem]{Example}
\newtheorem{rem}[Theorem]{Remark}
\newtheorem{setup}[Theorem]{Setup}

\theoremstyle{remark}
\newtheorem*{Acknowledgments}{Acknowledgments}

\newcommand{\depth}{\mathrm{depth}}

\newcommand{\Ker}{\mathrm{Ker}}

\def\cocoa{{\hbox{\rm C\kern-.13em o\kern-.07em C\kern-.13em o\kern-.15em A}}}

\usepackage{enumerate, comment}

\newcommand{\calF}{\mathcal{F}}

\newcommand{\calR}{\mathcal{R}}

\newcommand{\fka}{\mathfrak{a}}

\newcommand{\fkm}{\mathfrak{m}}
\newcommand{\fkn}{\mathfrak{n}}

\newcommand{\fkq}{\mathfrak{q}}

\def\depth{\operatorname{depth}}

\def\grade{\mathrm{grade}}

\def\Ker{\mathrm{Ker}}

\def\ol{\overline}

\def\red{\operatorname{red}}

\def\Gen{\operatorname{Gen}}
\def\MinGen{\operatorname{MinGen}}

\begin{document}

\title{Defining ideals of Cohen-Macaulay fiber cones}

\author{Reza Abdolmaleki}
\address{Department of Mathematics, Lahore University of Management Sciences, DHA, Lahore Cantt. 54792, Lahore, Pakistan}
\email{reza.abdolmaleki@lums.edu.pk}
\email{reza.abd110@gmail.com}
\author{Shinya Kumashiro}
\address{Department of Mathematics, Osaka Institute of Technology, 5-16-1 Omiya, asahi-ku, Osaka, 535-8585, Japan}
\email{shinya.kumashiro@oit.ac.jp}
\email{shinyakumashiro@gmail.com}

\begin{abstract} 
Let $A$ be a commutative Noetherian local ring with maximal ideal $\mathfrak{m}$, and let $I$ be an ideal. The fiber cone is then an image of the polynomial ring over the residue field $A/\mathfrak{m}$. The kernel of this map is called the defining ideal, and it is natural to ask how to compute it. In this paper, we provide a construction for the defining ideals of Cohen-Macaulay fiber cones.
\end{abstract}


\thanks{The second author was supported by JSPS KAKENHI Grant Number 24K16909.}

\maketitle

\section{Introduction}
The fiber cone of an ideal harbors significant importance in commutative algebra and algebraic geometry, offering insights from various perspectives.
Let $A$ be a commutative Noetherian local ring with the unique maximal ideal $\fkm$. We denote by $K$ the residue field $A/\fkm$. Let $A[t]$ be the polynomial ring in the variable $t$ over $A$. For an ideal $I$ of $A$, the \textit{Rees algebra} of $I$ is defined as
\[
\calR(I) := A[It],
\]
and the \textit{fiber cone} of \( I \) with respect to \( \fkm \) is defined as
\[
\calF(I) := \calR(I)/\fkm \calR(I).
\]
By definition, we have $\calF(I) \cong \bigoplus_{i \geq 0} I^i / \fkm I^i$. Thus, the Hilbert function of \(\mathcal{F}(I)\) provides insights into the minimal number of generators of the powers of \(I\). Moreover, the Krull dimension of \(\mathcal{F}(I)\), called the \textit{analytic spread} of \(I\), coincides with the minimal number of generators of any minimal reduction \(J\) of \(I\). 
From a geometric viewpoint, the role of $\mathcal{F}(I)$ in the process of blowing up $\text{Spec}(R)$ along $V (I)$ is pivotal. Indeed, $\text{Proj}(\mathcal{F}(I))$ corresponds to the fiber over the closed point of the blowup of $\text{Spec}(R)$ along $V (I)$. For further information about fiber cones, one can see, for example, \cite{HIO} and \cite{RV}. 

Setting \( I = (x_1, \dots, x_n) \) for \( x_1, \dots, x_n \in A \), we have \(\calF(I) = K[\overline{x_1}t, \dots, \overline{x_n} t]\), where \(\overline{x_i}t\) denotes the image of \( x_i t \in \calR(I) \) in \(\calF(I)\). We then consider the following canonical surjective ring homomorphism:
\[
\varphi: K[X_1, \dots, X_n] \to \calF(I), \quad X_i \mapsto \overline{x_i} t.
\]
The ideal \(\Ker \varphi\) is called the \textit{defining ideal} of the fiber cone \(\calF(I)\) of \( I \). Determining the defining ideal \(\Ker \varphi\) is crucial for understanding the structure of \(\calF(I)\), but it is often challenging. 
This problem is investigated for some certain classes of ideals $I$. For example, if $I$ is an ideal generated by $\dim A+1$ elements and $\grade (G_{+}(I))\ge \dim A-1$, then $\calF(I)$ is a hypersurface (\cite[Theorem 5.6]{HK}). If $I$ is a Freiman ideal, which is a special class of monomial ideals in the polynomial ring $A$, then the defining ideal of $\calF(I)$ has a $2$-linear free resolution (\cite[Theorem 1.3(g)]{HHZ}). In the papers \cite{HQS} and \cite{HZ}, the authors determine the defining ideals of fiber cones for special classes of monomial ideals $I$ in the polynomial ring with two variables (concave ideals, convex ideals, and certain symmetric ideals), see \cite[Theorems 2.6, 2.10, 3.5]{HQS} and \cite[Theorems 3.1, 3.2, 4.1, 4.2, and 4.3]{HZ}. One can find further results in \cite{AZ, HS, S}.






In this paper, we construct a certain ideal $\fka$ in $K[X_1, \dots, X_n]$, which is contained in $\Ker \varphi$ (Proposition \ref{lem34}). Then, we prove that $\fka = \Ker \varphi$ if $\calF(I)$ is Cohen-Macaulay (Theorem \ref{thm36}). We also give an example illustrating Theorem \ref{thm36} (Example \ref{exmonomial}).
Note that the Hilbert functions of Cohen-Macaulay fiber cones $\calF(I)$ are well-known (\cite[Theorem 6]{Shah1}). In the next section we prove the main result. In what follows, for a module $M$, $\ell(M)$ (resp, $\mu(M)$) denotes the length of $M$ (resp. the minimal number of generators of $M$). 

\if
we determine the defining ideals of Cohen-Macaulay fiber cones. Let $I$ be an $\fkm$-primary ideal of a Noetherian local ring $(A, \fkm, K)$ of dimension $d$, with the unique maximal ideal $\fkm$ and residue field $K$. We suppose the existence of a parameter ideal $Q$ of $A$ such that $Q$ is a reduction of $I$ (such an ideal always exists if the residue field is infinite). We set $Q=(x_1, \dots, x_d)$ and $I=Q+(x_{d+1}, \dots, x_{n})$, where $n=\mu(I)$ is the minimal number of generators of $I$. We then define the homomorphism $\varphi: K[X_1, \dots, X_n] \to \calF(I)$ by $X_i \mapsto \ol{x_i} t$. Let $r$ be the reduction number of $I$ with respect to $Q$ and derive the minimal generators of $I^i/QI^{i-1}$, denoted as $y_{i1}, y_{i2}, \dots, y_{iu_i}$, from $\Gen(I^i)$. This setup enables us (Theorem~\ref{thm36})to demonstrate that the kernel of $\varphi$ is generated by the polynomials 
\[
X_{\ell_1}X_{\ell_2} \cdots X_{\ell_{r+1}} - \sum_{\begin{smallmatrix}1 \le p \le d, \\ 1 \le q_1, \dots, q_{r} \le n \end{smallmatrix}} \overline{b_{p,q_1, \dots, q_{r}}} X_p X_{q_1}\cdots X_{q_{r}}.
\] 
for all $2\le i \le r+1$ and $ x_{\ell_1}x_{\ell_2} \cdots x_{\ell_i} \in \Gen(I^i)\setminus \MinGen(I^i/QI^{i-1})$.
\fi



\section{Main theorem}

In this section, we present a specific construction method for the defining ideal of Cohen-Macaulay fiber cones. Our journey begins with the following lemma, which is a key ingredient in our result.

\begin{lem}\label{lem21}
Let \( R = \bigoplus_{i \ge 0} R_i \) be a positively graded Noetherian ring such that \( (R_0, \fkn) \) is a local ring. 
Consider the exact sequence of finitely generated graded \( R \)-modules:
\[ 0 \to L \to M \xrightarrow{\varphi} N \to 0. \]
Suppose that \( N \) is a Cohen-Macaulay \( R \)-module with \( \dim_R N = \dim_R M \), and let \( \fkq \) be a graded parameter ideal of \( N \). Then, \(\varphi\) is an isomorphism if and only if \(\ell_R(M/\fkq M) = \ell_R(N/\fkq N)\).
\end{lem}

\begin{proof}
Using the snake lemma recursively, we obtain the exact sequence 
\[
0 \to L/\fkq L \to M/\fkq M \xrightarrow{\ol{\varphi}} N/\fkq N \to 0.
\]
Since all modules in the short exact sequence are of dimension zero, we can observe that $L/\fkq L=0 \iff \ell_R(L/\fkq L)=0 \iff \ell_R(M/\fkq M) = \ell_R(N/\fkq N)$. By Nakayama's Lemma, \( L = 0 \) if and only if \( L/\fkq L = 0 \); this implies the assertion.
\end{proof}

Here, we illustrate our strategy for the main theorem. Consider a surjective graded ring homomorphism \(\varphi: K[X_1, \dots, X_n] \to \calF(I)\) defined by \(X_i \mapsto \ol{x_i} t\). We aim to find an ideal \(\fka \subseteq \Ker \varphi\) as a candidate for \(\Ker \varphi\). This leads to the induced exact sequence
\[
0 \to \Ker\varphi/ \fka \to K[X_1, \dots, X_n]/\fka \xrightarrow{\ol{\varphi}} \calF(I) \to 0
\]
as \(K[X_1, \dots, X_n]/\fka\)-modules. Suppose that $\fka$ satisfies \(\dim K[X_1, \dots, X_n]/\fka = \dim \calF(I) = d\), and assume that \(\calF(I)\) is Cohen-Macaulay. Then, by Lemma \ref{lem21}, we can determine if \(\varphi\) is an isomorphism by checking that the equality \(\ell_R(\calF(I)/ \fkq \calF(I)) = \ell_R(K[X_1, \dots, X_n]/(\fka+\fkq))\) holds, where \(\fkq\) is a graded parameter ideal of \(\calF(I)\).

To illustrate the ideal \(\fka\) mentioned above, we prepare the following notations.

\begin{setup}\label{setup22}
In what follows, let $(A, \fkm, K)$ be a Noetherian local ring of dimension $d$ with the unique maximal ideal $\fkm$ and the residue field $K$. Let $I$ be an $\fkm$-primary ideal of $A$. 
We assume that there exists a parameter ideal $Q$ of $A$ such that $Q$ is a reduction of $I$, that is, $I^{m+1}=QI^m$ for some $m>0$. Note that if $K$ is infinite, such a parameter ideal $Q$ always exists (\cite[Proposition 4.6.8]{BH}). We write $Q=(x_1, \dots, x_d)$. Since $x_1, \dots, x_d$ is a part of minimal generators of $I$, we can write $I=Q+(x_{d+1}, \dots, x_{n})$, where $n=\mu(I)$ is the minimal number of generators of $I$. We consider the defining ideal obtained by 
\[
\varphi: K[X_1, \dots, X_n] \to \calF(I); X_i \mapsto \ol{x_i} t.
\]

We set 
\[
r:=\red_Q(I)=\min\{m \ge 0 : I^{m+1} = QI^m \},
\]
the {\it reduction number} of $I$ with respect to $Q$. Let $u_i=\mu(I^i/QI^{i-1})$ for $i\ge 1$. 
For $i\ge 1$, we denote by $\Gen(I^i)$ the set of (not necessarily minimal) generators of $I^i$:
\[
\Gen(I^i):=\{x_{\ell_1}x_{\ell_2} \cdots x_{\ell_i} : 1 \le \ell_s \le n \text{ for $1 \le s \le i$}\}.
\]
We then choose a system of minimal generators $y_{i1}, y_{i2}, \dots, y_{iu_i}$ of $I^i/QI^{i-1}$ from $\Gen(I^i)$, and set 
\[
\MinGen(I^i/QI^{i-1}):=\{y_{i1}, y_{i2}, \dots, y_{iu_i}\}\subseteq \Gen(I^i)
\] 
for $i\ge 2$. 

\end{setup}

\begin{rem}\label{rem31}
By the definition of $r$, $u_i=0$ and $\MinGen(I^i/QI^{i-1})$ is the empty set if $i\ge r+1$. On the other hand, if $2 \le i \le r$, the choice of $y_{i1}, y_{i2}, \dots, y_{iu_i}$ is not necessarily unique. 
\end{rem}


With the Setup \ref{setup22} we construct the generators of the defining ideal. Assume that $2 \le i \le r+1$. 
For each element $x_{\ell_1}x_{\ell_2} \cdots x_{\ell_i} \in \Gen(I^i)\setminus \MinGen(I^i/QI^{i-1})$, there exists $a_j\in A$ such that 
\begin{align}\label{eqadd}
x_{\ell_1}x_{\ell_2} \cdots x_{\ell_i} - \sum_{j=1}^{u_i} a_j y_{ij} \in QI^{i-1}
\end{align}
since $x_{\ell_1}x_{\ell_2} \cdots x_{\ell_i} \in (y_{i1}, y_{i2}, \dots, y_{iu_i}) +QI^{i-1}$ by the definition of $y_{ij}$. 

\begin{rem}\label{rem32}
In the case where $i=r+1$, we choose $a_j =0$ for all $1\le j \le u_i$. This is possible since $x_{\ell_1}x_{\ell_2} \cdots x_{\ell_{r+1}} \in I^{r+1} =QI^r$. 
\end{rem}

Since $QI^{i-1}$ is generated by elements $x_p x_{q_1}\cdots x_{q_{i-1}}$ for all $1 \le p \le d$ and $1 \le q_1, \dots, q_{i-1} \le n$, by \eqref{eqadd}, we can write 
\[
x_{\ell_1}x_{\ell_2} \cdots x_{\ell_i} - \sum_{j=1}^{u_i} a_j y_{ij} = \sum_{\begin{smallmatrix}1 \le p \le d, \\
1 \le q_1, \dots, q_{i-1} \le n \end{smallmatrix}} b_{p,q_1, \dots, q_{i-1}} x_p x_{q_1}\cdots x_{q_{i-1}}
\]
for some $b_{p,q_1, \dots, q_{i-1}} \in A$. Letting $\ol{a}$ denote the image of $a\in A$ into $K=A/\fkm$ and $Y_{ij}$ denote a monomial of degree $i$ in $K[X_1, \dots, X_n]$ with $\varphi(Y_{ij}) = \ol{y_{ij}} t^{i}$ (recalling that $y_{ij}\in \Gen(I^i)$ and $\varphi(X_i) = \ol{x_i}t$, such a monomial always exists), we then define
\begin{align}\label{eqz}
Z_{\ell_1, \dots, \ell_i} := X_{\ell_1}X_{\ell_2} \cdots X_{\ell_i} - \sum_{j=1}^{u_i} \ol{a_j} Y_{ij} - \sum_{\begin{smallmatrix}1 \le p \le d, \\
1 \le q_1, \dots, q_{i-1} \le n \end{smallmatrix}} \ol{b_{p,q_1, \dots, q_{i-1}}} X_p X_{q_1}\cdots X_{q_{i-1}}
\end{align}
in $K[X_1, \dots, X_n]$ for each $x_{\ell_1}x_{\ell_2} \cdots x_{\ell_i} \in \Gen(I^i)\setminus \MinGen(I^i/QI^{i-1})$. 

\begin{rem}\label{rem33}
$Z_{\ell_1, \dots, \ell_i}$ is a homogeneous element with $\deg Z_{\ell_1, \dots, \ell_i} =i$.
\end{rem}

\begin{prop}\label{lem34}
For each $2\le i \le r+1$ and $x_{\ell_1}x_{\ell_2} \cdots x_{\ell_i} \in \Gen(I^i)\setminus \MinGen(I^i/QI^{i-1})$, $Z_{\ell_1, \dots, \ell_i} \in \Ker \varphi$. 
\end{prop}

\begin{proof}
This follows by the construction of $Z_{\ell_1, \dots, \ell_i}$: 
\begin{align*}
\varphi(Z_{\ell_1, \dots, \ell_i}) =& \varphi(X_{\ell_1}X_{\ell_2} \cdots X_{\ell_i} - \sum_{j=1}^{u_i} \ol{a_j} Y_{ij} - \sum_{\begin{smallmatrix}1 \le p \le d, \\
1 \le q_1, \dots, q_{i-1} \le n \end{smallmatrix}} \ol{b_{p,q_1, \dots, q_{i-1}}} X_p X_{q_1}\cdots X_{q_{i-1}})\\
=& \left(\ol{x_{\ell_1}x_{\ell_2} \cdots x_{\ell_i}} - \sum_{j=1}^{u_i} \ol{a_j y_{ij}} - \sum_{\begin{smallmatrix}1 \le p \le d, \\
1 \le q_1, \dots, q_{i-1} \le n \end{smallmatrix}} \ol{b_{p,q_1, \dots, q_{i-1}} x_p x_{q_1}\cdots x_{q_{i-1}}}\right)t^i \\
=&0
\end{align*}
\end{proof}

Set 
\[
\fka := (Z_{\ell_1, \dots, \ell_i} : 2\le i \le r+1, x_{\ell_1}x_{\ell_2} \cdots x_{\ell_i} \in \Gen(I^i)\setminus \MinGen(I^i/QI^{i-1}))
\]
as an ideal of $K[X_1, \dots, X_n]$. By Proposition \ref{lem34}, we obtain that $\fka\subseteq \Ker\varphi$. 

With the above notations, we get the following proposition. 

\begin{prop}\label{prop35}
The following statements hold true:
\begin{enumerate}[\rm(i)]
\item $\fka$ is a graded ideal generated by elements of degree $2 \le i \le r+1$.
\item $\fka + (X_1, \dots, X_d)$ contains $(X_1, \dots, X_n)^{r+1}$.
\item $\dim K[X_1, \dots, X_n]/\fka = d$ and $(X_1, \dots, X_d)$ is a parameter ideal of $K[X_1, \dots, X_n]/\fka$.
\end{enumerate}
\end{prop}

\begin{proof}
(i): This follows from Remark \ref{rem33}. 

(ii): Let $i=r+1$. Since $\MinGen(I^{r+1}/QI^{r})$ is empty (Remark \ref{rem31}), we have $Z_{\ell_1, \dots, \ell_{r+1}}$ for each $1 \le \ell_1, \dots, \ell_{r+1} \le n$. 
By Remark \ref{rem32}, we have 
\[
Z_{\ell_1, \dots, \ell_{r+1}} = X_{\ell_1}X_{\ell_2} \cdots X_{\ell_{r+1}} - \sum_{\begin{smallmatrix}1 \le p \le d, \\
1 \le q_1, \dots, q_{r} \le n \end{smallmatrix}} \ol{b_{p,q_1, \dots, q_{r}}} X_p X_{q_1}\cdots X_{q_{r}}.
\]
Since $1 \le p \le d$, it follows that 
\[
X_{\ell_1}X_{\ell_2} \cdots X_{\ell_{r+1}} = Z_{\ell_1, \dots, \ell_{r+1}} + \sum_{\begin{smallmatrix}1 \le p \le d, \\
1 \le q_1, \dots, q_{r} \le n \end{smallmatrix}} \ol{b_{p,q_1, \dots, q_{r}}} X_p X_{q_1}\cdots X_{q_{r}} \in \fka +(X_1, \dots, X_d). 
\]
Hence, all monomials of degree $r+1$ are in $\fka +(X_1, \dots, X_d)$, that is, $\fka + (X_1, \dots, X_d)$ contains $(X_1, \dots, X_n)^{r+1}$. 

(iii): Since we have a surjection $K[X_1, \dots, X_n]/\fka \xrightarrow{\ol{\varphi}} \calF(I)$ and $\dim \calF(I) =d$ (see, for example, \cite[Exercise 4.6.13(b)]{BH}), $\dim K[X_1, \dots, X_n]/\fka \ge \dim \calF(I) =d$. On the other hand, since $\dim K[X_1, \dots, X_n]/(\fka + (X_1, \dots, X_d)) =0$ by (ii), the Krull's height theorem shows that $\dim K[X_1, \dots, X_n]/\fka \le d$. Therefore, we obtain that  $\dim K[X_1, \dots, X_n]/\fka =d$. Thus, $(X_1, \dots, X_d)$ is a parameter ideal of $K[X_1, \dots, X_n]/\fka$.
\end{proof}

\begin{thm}\label{thm36}
If $\calF(I)$ is a Cohen-Macaulay ring, then $\Ker \varphi =\fka$. 
\end{thm}

\begin{proof} 
We have $\fka\subseteq \Ker\varphi$ by Proposition \ref{lem34}. This leads to the exact sequence
\[
0 \to \Ker\varphi/ \fka \to K[X_1, \dots, X_n]/\fka \xrightarrow{\ol{\varphi}} \calF(I) \to 0
\]
as \(K[X_1, \dots, X_n]/\fka\)-modules. On the other hand, \(\dim K[X_1, \dots, X_n]/\fka = d= \dim \calF(I)\), and \((X_1, \dots, X_d)\calF(I) = (\ol{x_1}t, \ol{x_2}t, \dots, \ol{x_d}t)\) is a parameter ideal of \(\calF(I)\) by Proposition \ref{prop35}(iii). Hence, since \(\calF(I)\) is Cohen-Macaulay, by Lemma \ref{lem21}, we only need to prove that the following equation holds:
\begin{align}\label{eqlength}
\ell(K[X_1, \dots, X_n]/(\fka +(X_1, \dots, X_d))) = \ell(\calF(I)/(\ol{x_1}t, \ol{x_2}t, \dots, \ol{x_d}t)).
\end{align}

Indeed, we have 
\begin{align}
\begin{split}\label{rhs}
(\text{right-hand side of (\ref{eqlength})}) &= \ell(\calF(I)/Qt\calF(I)) = \ell\left(\bigoplus_{i\ge 0} I^i/(\fkm I^i + QI^{i-1})\right) \\
&= \sum_{i\ge 0} \ell(I^i/(\fkm I^i + QI^{i-1})) \\
&= 1 + \sum_{i\ge 1} u_i \\
&= 1 + \sum_{i=1}^{r} u_i,
\end{split}
\end{align}
where the last equation follows from Remark \ref{rem31}.

On the other hand, denoting by $[R]_i$ the $i$th component of a graded ring $R$, we also have
\begin{align}
\begin{split}\label{lhs}
\text{(left-hand side of \eqref{eqlength})} &= \ell\left(\bigoplus_{i \ge 0} [K[X_1, \dots, X_n]/(\fka + (X_1, \dots, X_d))]_i\right) \\
&= \sum_{i \ge 0} \ell\left([K[X_1, \dots, X_n]/(\fka + (X_1, \dots, X_d))]_i\right) \\
&= 1 + (n-d) + \sum_{i \ge 2} \ell\left([K[X_1, \dots, X_n]/(\fka + (X_1, \dots, X_d))]_i\right) \\
&= 1 + (n-d) + \sum_{i=2}^{r} \ell\left([K[X_1, \dots, X_n]/(\fka + (X_1, \dots, X_d))]_i\right),
\end{split}
\end{align}
where the third equation follows from Proposition \ref{prop35}(i) and the fourth equation follows from Proposition \ref{prop35}(ii).

Let \(2 \le i \le r\) and \(1 \le \ell_1, \dots, \ell_i \le n\). Suppose that \(X_{\ell_1} X_{\ell_2} \cdots X_{\ell_i}\) is a monomial of degree \(i\) such that \(x_{\ell_1} x_{\ell_2} \cdots x_{\ell_i} \not\in \MinGen (I^i/QI^{i-1})\). Recall that in this case, we can define \(Z_{\ell_1, \dots, \ell_i} \in \fka\) (see \eqref{eqz}). Thus, we obtain that
\begin{align*}
X_{\ell_1}X_{\ell_2} \cdots X_{\ell_i} - \sum_{j=1}^{u_i} \ol{a_j} Y_{ij} 
=& Z_{\ell_1, \dots, \ell_i} + \sum_{\begin{smallmatrix}1 \le p \le d, \\
1 \le q_1, \dots, q_{i-1} \le n \end{smallmatrix}} \ol{b_{p,q_1, \dots, q_{i-1}}} X_p X_{q_1} \cdots X_{q_{i-1}} \in \fka + (X_1, \dots, X_d).
\end{align*}

It follows that \(X_{\ell_1} X_{\ell_2} \cdots X_{\ell_i} \in \fka + (X_1, \dots, X_d) + (Y_{i1}, \dots, Y_{iu_i})\), that is, \(Y_{i1}, \dots, Y_{iu_i}\) span \([K[X_1, \dots, X_n]/(\fka + (X_1, \dots, X_d))]_i\) as a \(K\)-vector space. This concludes that \(\ell([K[X_1, \dots, X_n]/(\fka + (X_1, \dots, X_d))]_i) \le u_i\) for all \(2 \le i \le r\). Hence, we obtain that
\begin{align}
\begin{split}\label{eqfinal}
\text{(left-hand side of \eqref{eqlength})} &= 1 + (n-d) + \sum_{i=2}^{r} \ell\left([K[X_1, \dots, X_n]/(\fka + (X_1, \dots, X_d))]_i\right) \\
&= 1 + u_1 + \sum_{i=2}^{r} \ell\left([K[X_1, \dots, X_n]/(\fka + (X_1, \dots, X_d))]_i\right) \\
&\le 1 + u_1 + \sum_{i=2}^{r} u_i \\
&= \text{(right-hand side of \eqref{eqlength})},
\end{split}
\end{align}
where the first equality follows from \eqref{lhs}, the third inequality follows from the inequalities 
\[
\ell([K[X_1, \dots, X_n]/(\fka + (X_1, \dots, X_d))]_i) \le u_i
\]
 for all \(2 \le i \le r\), and the last equality follows from \eqref{rhs}. Recalling that we have a surjection \(\ol{\varphi}: K[X_1, \dots, X_n]/(\fka + (X_1, \dots, X_d)) \to \calF(I)/(\ol{x_1}t, \ol{x_2}t, \dots, \ol{x_d}t)\), the left-hand side of \eqref{eqlength} $\ge$ the right-hand side of \eqref{eqlength}. It follows that the inequality \eqref{eqfinal} must be an equality. Thus, we complete the proof.
\end{proof}

We illustrate our theorem with the following example. 

\begin{ex}\label{exmonomial}
Let $A=K[[x,y]]$ be the formal power series ring over a field $K$ with maximal ideal $\fkm=(x,y)$. Set 
\begin{center}
$I=(x^7, x^5y, x^4y^2, x^2y^6, y^{12})$ \quad and \quad $Q=(x^7+x^4y^2+y^{12}, x^5y+x^2y^6)$. 
\end{center} 
For simplicity, we write $f=x^7+x^4y^2+y^{12}$ and $g= x^5y+x^2y^6$. We then consider the ring homomorphism $\varphi: K[X_1, \dots, X_5] \to \calF(I)$ defined by 
\begin{center}
$\varphi(X_1) = \ol{f}t$, \ \ $\varphi(X_2) = \ol{g}t$, \ \ $\varphi(X_3) = \ol{x^7}t$, \ \ $\varphi(X_4) = \ol{x^2y^6}t$, \ \ $\varphi(X_5) = \ol{y^{12}}t$.
\end{center} 
Then, the following hold true. 
\begin{enumerate}[\rm(i)] 
\item $Q$ is a reduction of $I$ with $I^2=QI+(x^4y^{12})$ and $I^3=QI^2$. 
\item $\fkm I^2=\fkm QI$ and $\calF(I)$ is a Cohen-Macaulay ring. 
\item $\Ker \varphi =(X_3^2-X_1X_3, X_3X_4, X_3X_5, X_4X_5-X_2X_5, X_5^2-X_1X_5, X_4^3+X_2^2X_4-2X_2X_4^2)$. 
\end{enumerate}
\end{ex}

\begin{proof}
(i): Set $J= (x^7, x^2y^6, y^{12})$. Since $I=Q + J$, we have $I^2 = Q^2 + QJ + J^2 = QI+J^2$. For generators $x^{14}, x^9y^6,  x^4y^{12}, x^2y^{18}, y^{24}$ of $J^2$ (not including $x^7y^{12}$), we can observe that 
\begin{align}\label{eq291}
\begin{split}
x^{14} =& x^7f-x(x^5y)^2 - y^5 (x^5y)(x^2y^6)\in QI+\fkm I^2, \\
x^9y^6 =& y^3 (x^5y)(x^4y^2)\in \fkm I^2, \\
x^7y^{12} =& y^5(x^5y)(x^2y^6) \in \fkm I^2,\\
x^2y^{18} =& gy^{12} -xy(x^2y^6)^2 \in QI +\fkm I^2, \text{ and}\\
y^{24} =& fy^{12} - y^5(x^5y)(x^2y^6) - x^2(x^2y^6)^2 \in QI + \fkm I^2.
\end{split}
\end{align}
Hence, $I^2=QI+J^2 \subseteq QI+ (x^7y^{12}) + \fkm I^2 \subseteq I^2$. Using Nakayama's lemma, we get $I^2=QI+(x^4y^{12})$. Multiplying both sides of the equality by $I$, we obtain that 
\[
I^3= QI^2 + x^4y^{12} I = QI^2 + x^4y^{12} (Q+J) = QI^2 + x^4y^{12} J.
\]
For generators $x^{11}y^{12}, x^6y^{18}, x^{4}y^{24}$ of $x^4y^{12} J$, we can observe that 
\begin{align}\label{eq292}
\begin{split}
x^{11}y^{12} =&y^3(x^2y^6)(x^4y^2)(x^5y) \in \fkm I^3, \\
x^6y^{18} =& -g^2(x^2y^6) + 2g(x^2y^6)^2 + y^2(x^4y^2)^3 \in QI^2 +\fkm I^3, \text{ and}\\
x^{4}y^{24} =& g(x^2y^6)(y^{12}) -xy (x^2y^6)^3 \in QI^2+\fkm I^3.
\end{split}
\end{align}
Hence, $I^3=QI^2 + x^4y^{12} J \subseteq QI^2+\fkm I^3 \subseteq I^3$. Now, Nakayama's lemma implies that $I^3=QI^2$.

(ii): By (i), we have $\fkm I^2 = \fkm QI +x^4y^{12}\fkm$. For generators $x^5y^{12}, x^4y^{13}$ of $x^4y^{12}\fkm$, we can observe that 
\begin{align*}
x^5y^{12} =& xg(x^2y^6) - y^3 (x^4y^2)^2 \in \fkm QI +\fkm^2 I^2 \text{ and}\\
x^4y^{13} =& (yg-xf)(x^2y^6) + x^2y^2(x^4y^2)^2 +xg(y^{12}) -y^5(x^4y^2)(x^2y^6) \in \fkm QI +\fkm^2 I^2. 
\end{align*}
Hence, $\fkm I^2 = \fkm QI +x^4y^{12}\fkm \subseteq \fkm QI +\fkm^2 I^2 \subseteq \fkm I^2$. By Nakayama's lemma, we have $\fkm I^2=\fkm QI$. 
In particular, $I^2/QI$ is a $K$-vector space, and this is of dimension one by (i). That is, $\ell(I^2/QI) =1$. This shows that $I$ is a Sally ideal in the sense of \cite{JPV} (or $I$ has {\it almost minimal multiplicity} in the sense of \cite{RV}). Therefore, by \cite[Theorem 3.4]{JPV}, $\calF(I)$ is Cohen-Macaulay. 

(iii): By the first equation in \eqref{eq291}, we have $(x^{7})^2 - x^7f = -x(x^5y)^2  - y^5 (x^5y)(x^2y^6)\in \fkm I^2$. It follows that $X_3^2 -X_3X_1 \in \Ker \varphi$. In the same way, we observe that 
\[
Z_{33}:= X_3^2 -X_3X_1, \quad Z_{34}:= X_3X_4, \quad Z_{35}:=X_3X_5, \quad Z_{45}:= X_4X_5 - X_2X_5, \quad Z_{55}:= X_5^2-X_1X_5
\] 
are in $\Ker \varphi$. Similarly, by \eqref{eq292}, we obtain that $Z_{444}:= X_4^3+X_2^2X_4-2X_2X_4^2$ are in $\Ker \varphi$. 

Set $\fka:=(Z_{33}, Z_{34}, Z_{35}, Z_{45}, Z_{55}, Z_{444})$, and consider the exact sequence 
\[
0 \to \Ker\varphi/ \fka \to K[X_1, \dots, X_5]/\fka \xrightarrow{\ol{\varphi}} \calF(I) \to 0
\]
as $K[X_1, \dots, X_n]/\fka$-modules. Since $\calF(I)$ is Cohen-Macaulay by (ii), and $\varphi(X_1) = \ol{f}t, \varphi(X_2) = \ol{g}t$ is a regular sequence of $\calF(I)$, we obtain that 
\[
0 \to \Ker\varphi/ (\fka +(X_1, X_2)) \to K[X_1, \dots, X_5]/(\fka+(X_1, X_2)) \to \calF(I)/(\ol{f}t, \ol{g}t) \to 0.
\]
It is straightforward to check that $\fka+(X_1, X_2) = (X_1, X_2) + (X_3^2, X_3X_4, X_3X_5, X_4X_5, X_5^2, X_4^3)$; hence, $\ell(K[X_1, \dots, X_5]/(\fka+(X_1, X_2))) = 5$. On the other hand, we get $\ell(\calF(I)/(\ol{f}t, \ol{g}t)) = \sum_{i\ge 0} \mu(I^i/QI^{i-1}) = 1 + \mu(I/Q) + \mu(I^2/QI) = 5$ by (i). Therefore, $\Ker\varphi/ (\fka +(X_1, X_2)) =0$. By the graded Nakayama's lemma, this proves that $\fka = \Ker\varphi$. 
\end{proof}

\begin{rem}
We note that in Example \ref{exmonomial}(i), the equality $I^3=QI^2$ does not hold in the polynomial ring $K[x,y]$, although $I^3\subseteq QI^2+\fkm I^3$ holds. One can check this by CoCoA \cite{Co}. 
\end{rem}

The following example demonstrates that the ideal $\fka$ in Theorem \ref{thm36} does not necessarily equal $\Ker\varphi$ when $\calF(I)$ is not a Cohen-Macaulay ring.

\begin{ex}
Let $A=K[[x,y]]$ be the formal power series ring over a field $K$ with maximal ideal $\fkm=(x,y)$. Set 
\begin{center}
$I=(x^5, x^4y, x^3y^3, xy^4, y^5)$ \quad and \quad $Q=(x^5, y^5)$. 
\end{center} 
We then consider the ring homomorphism $\varphi: K[X_1, \dots, X_5] \to \calF(I)$ defined by 
\begin{center}
$\varphi(X_1) = \ol{x^5}t$, \ \ $\varphi(X_2) = \ol{x^4y}t$, \ \ $\varphi(X_3) = \ol{x^3y^3}t$, \ \ $\varphi(X_4) = \ol{xy^4}t$, \ \ $\varphi(X_5) = \ol{y^{5}}t$.
\end{center} 
Then, the following hold true. 
\begin{enumerate}[\rm(i)] 
\item $Q$ is a reduction of $I$ with $I^2=QI+(x^8y^2,x^2y^8)$, $I^3=QI^2 + (x^{12}y^3,x^3y^{12})$, and $I^4=QI^3$. 
\item $\fka \subsetneq \Ker \varphi$. Indeed, $\fka = (X_2X_3, X_2X_4 - X_1X_5, X_3^2, X_3X_4, X_2^4-X_1^3X_4, X_4^4-X_2X_5^3)$, but \\
$\Ker \varphi = \fka +(X_1X_3, X_3X_5, X_1^2X_4^2 - X_2^3X_5, X_1X_4^3-X_2^2X_5^2)$. 
\item $\depth \calF(I) =0$. 
\end{enumerate}
\end{ex}

\begin{proof}
(i): It is straightforward to check. 

(ii): We observe that 
\begin{align*}
(x^4y)(x^3y^3) - x(x^5)(xy^4)=0, \quad (x^4y)(xy^4) - (x^5)(y^5) =0, \quad (x^3y^3)^2 -xy(x^5)(y^5)=0, \\
(x^3y^3)(xy^4) - y(x^4y)(y^5) =0, \quad (x^4y)^4-(x^5)^3(xy^4) =0, \quad \text{and} \quad (xy^4)^4-(x^4y)(y^5)^3 =0.
\end{align*}
Hence, $\fka$ is what we desired. On the other hand, one can compute by CoCoA \cite{Co} the generators of $\Ker \varphi$ is as stated in (ii). 

(iii): Since the generators of the defining ideal are computed, one can check again by CoCoA that $\depth \calF(I) =0$. 
\end{proof}

\begin{Acknowledgments}
The authors first came together in 2019 in Essen, where they had the opportunity to meet J\"{u}rgen Herzog. They wish to pay tribute to Herzog for the invaluable knowledge and inspiration he imparted to them. They are deeply grateful for the privilege of learning from him.
\end{Acknowledgments}



\begin{thebibliography}{99}
\bibitem{Co} 
{\sc J. Abbott, A. M. Bigatti, L. Robbiano}, CoCoA: a system for doing Computations in Commutative Algebra. Available at: \url{http://cocoa.dima.unige.it}.


\bibitem{AZ} 
{\sc R. Abdolmaleki, R. Zaare-Nahandi}, Toric ideals which are determinantal, available at arXiv:2106.03900 


\bibitem{BH} 
{\sc W. Bruns, J. Herzog}, Cohen-Macaulay rings. Cambridge Studies in Advanced Mathematics, 39. Cambridge University Press, Cambridge, 1993. 


\bibitem{HK}
{\sc W. J. Heinzer, M.-K. Kim}, Properties of the fiber cone of ideals in local rings. {\em Comm. Algebra} {\bf 31} (2003), no. 7, 3529--3546.

\bibitem{HIO} 
{\sc M. Herrmann, S. Ikeda and U. Orbanz},  Equimultiplicity and Blowing up. Springer-Verlag, 1988.

\bibitem{HHZ} 
{\sc J. Herzog, T. Hibi, G. Zhu}. The relevance of Freiman's theorem for combinatorial commutative algebra,  {\em Math. Z.},  {\bf 291} no. 3-4, (2019), 999--1014.




\bibitem{HQS}  
{\sc J. Herzog, A. A. Qureshi,  M. M. Saem},  The fiber cone of a monomial ideal in two variables, {\em J. Symbolic Comput.} {\bf 94} (2019), 52--69. 



\bibitem{HZ}  
{\sc J. Herzog,  G. Zhu}, On the fiber cone of monomial ideals,  {\em Arch. Math.} {\bf 113} (2019),  469--481.


\bibitem{HS} 
{\sc R. H\"{u}bl, I. Swanson}, Normal cones of monomial primes,  {\em Math. Comp.}, {\bf 72} (2003), 459--475.





 

\bibitem{JPV}
{\sc A. V. Jayanthan, T. J. Puthenpurakal, J. K. Verma}, On fiber cones of $\fkm$-primary ideals, {\em Canad. J. Math.}, {\bf 59} (2007), no. 1, 109--126.



\bibitem{RV}
{\sc M. E. Rossi, G. Valla}, Hilbert functions of filtered modules. Lecture Notes of the Unione Matematica Italiana, 9. Springer-Verlag, Berlin; UMI, Bologna, 2010. 


\bibitem{Shah1}
{\sc K. Shah}, On the Cohen-Macaulayness of the fiber cone of an ideal. {\em J. Algebra}, {\bf 143} (1991), no. 1, 156--172.

\bibitem{Shah2}
{\sc K. Shah, On equimultiple ideals}, {\em Math. Z.}, {\bf 215} (1994), 13--24.

\bibitem{S} 
{\sc B. Sturmfels},   Gr\"obner Bases and Convex Polytopes.  Amer. Math. Soc., Providence, RI, 1995.

\end{thebibliography}
\end{document}